\documentclass[12pt, twoside, leqno]{article}

\usepackage{amsmath,amsthm}
\usepackage{amssymb}

\usepackage{enumerate}

\usepackage{graphicx}

\usepackage{pgf,tikz}
\usetikzlibrary{arrows}

\newtheorem{thm}{Theorem}[section]

\newtheorem{prop}[thm]{Proposition}

\theoremstyle{definition}

\newtheorem{rem}[thm]{Remark}

\numberwithin{equation}{section}

\usepackage{pgf,tikz}
\usetikzlibrary{arrows}

\begin{document}

\baselineskip=12.1pt


\title{ The line graph of the crown graph is distance integral}

\author{S. Morteza Mirafzal\\
Department of Mathematics \\
  Lorestan University, Khorramabad, Iran\\
E-mail: mirafzal.m@lu.ac.ir\\
E-mail: smortezamirafzal@yahoo.com}

\date{}

\maketitle

\renewcommand{\thefootnote}{}


\footnote{\emph{Date}:  }

\renewcommand{\thefootnote}{\arabic{footnote}}
\setcounter{footnote}{0}
\date{}

\begin{abstract} 
 The distance eigenvalues of a connected graph $G$ are the eigenvalues of its distance matrix
$D(G)$.  A graph is called distance integral if all of its
distance eigenvalues are integers. Let $n \geq 3$ be an integer.  A crown graph $Cr(n)$ is a graph obtained from the complete bipartite graph $K_{n,n}$ by removing a perfect matching. Let $L(Cr(n))$ denote the line graph of the crown graph $Cr(n)$.    In this paper,  by using the   orbit partition method in algebraic graph theory,   we determine the set of  all distance eigenvalues of $L(Cr(n))$ and show that this graph   is  distance  integral.
\end{abstract}

\maketitle

\section{ Introduction and Preliminaries}
  In this paper, a graph $G=(V,E)$ is
considered as an undirected simple graph where $V=V(G)$ is the vertex-set
and $E=E(G)$ is the edge-set. For all the terminology and notation
not defined here, we follow [3,4,5,6,7].\

Let $G=(V,E)$ be a graph and  $A=A(G)$ be an adjacency matrix of $G$. The characteristic polynomial of $G$ is defined as $P(G; x)=P(x) = |xI-A|$. A zero of $p(x)$ is called an eigenvalue of the graph $G$. A graph
is called integral, if all the eigenvalues are integers. The study of integral graphs  was initiated by Harary
and Schwenk in 1974 (see [8]). A survey of papers up to 2002 has been  appeared in [2], but
more than a hundred new studies on integral graphs have been published in the last
twenty  years. \newline 
Let $n$ be the number of vertices of the graph $G$. The distance matrix $D=D(G)$ is an $n \times n $ matrix indexed by $V$, such that $D_{u,v} = d_{G}(u, v)=d(u,v)$, where $d_{G}(u, v)$ is the distance between the vertices $u$ and $v$ in the graph $G$.
The characteristic polynomial $P(D; x) = |xI - D|=D_G(x)$ is the distance characteristic polynomial of $G$. Since $D$ is a real
symmetric matrix, the distance characteristic polynomial $D_G(x)$ has real zeros. Every zero of the polynomial $D_G(x)$ is called a distance eigenvalue of the graph $G$. 
A survey on the distance
spectra of graphs has been appeared in [1]. 
A graph $G$ is distance integral (briefly, $D$-integral) if all the distance eigenvalues  of $G$ are integers. Although there are many
 papers that study distance spectrum of graphs and their applications, the $D$-integral graphs are studied only in a few number of papers (see [14]). \
  
Let $n \geq 3$ be an integer.  A crown graph $Cr(n)$ is a graph obtained from the complete bipartite graph $K_{n,n}$ by removing a perfect matching.  The bipartite Kneser graph $H(n,k)$, $1 \leq k \leq n-1$, is a bipartite graph with the vertex-set consisting of all $k$-subsets and $(n-k)$-subsets of the set $[n]=\{1,2,3,\dots,n\}$, in which two vertices $v$ and $w$ are adjacent if and only if $v \subset w$ or $w \subset v$. Recently, this class of graphs has been studied form several  aspects [9,11,12,13]. 
It is easy to see that the crown graph $Cr(n)$ is isomorphic with the bipartite graph 
$H(n,1)$. The crown graph $Cr(n)$ is a vertex and edge-transitive graph of order $2n$ and regularity $n-1$ with diameter 3.\

 In [14] it has been shown that the graph $Cr(n)$ is a distance integral graph. 
 Let $L(Cr(n))$ denote the line graph of the crown graph $Cr(n)$.  It is not hard to see that the graph $L(Cr(n))$ is a vertex-transitive graph of order $n(n-1)$ and regularity $2(n-1)-2=2n-4$ with diameter 3.\

From the various interesting properties  of the graph $L(Cr(n))$, we interested in its distance eigenvalues.  In this paper, we wish to show that the graph $L(Cr(n))$  is distance integral. We explicitly determine all the distinct distance eigenvalues of the graph $L(Cr(n))$. The method which we use in this paper is completely different from what has been used in [14].  The main tool which we use in our work is the orbit partition method in
algebraic graph theory which we have already employed it in determining the adjacency eigenvalues of a particular family of  graph [10]. In this paper,  We show how we can find, by using this method, the set
of all distinct distance eigenvalues of the graph  $L(Cr(n))$.

The set of all permutations of a set $V$ is denoted by $Sym(V)$. A permutation group on $V$ is a subgroup of $Sym(V)$. If $G=(V,E) $ is a graph, 
 then we can view each automorphism as a
permutation of
$ V,$  and so $Aut (G)$ is a permutation group.
A permutation representation of a group $\Gamma$ is a homomorphism from $\Gamma$
into $Sym(V)$ for some set  $V.$ A permutation representation is also referred
to as an $action$  of $\Gamma $ on the set $V, $ in which case we say that $\Gamma$ acts on $V.$ A permutation group $\Gamma$ on  $V$  is $transitive$ if given any two elements $x$ and
$y$ from $V$ there is an element $g \in \Gamma$ such that $x^g=y.$  For each $v\in V,$  the set
$v^{\Gamma}=\{ v^g \ | \ g \in \Gamma  \}$ is called an $orbit$  of  $\Gamma$.  
 It is easy to see that if $\Gamma$ acts on $V$, then $\Gamma$ is
transitive  on $V$ (or $\Gamma$ acts $transitively$ on $V$), when there is just
one orbit. It is easy to see that the set of orbits of $\Gamma$ on $V$ is a partition of the set $V$.\\
 A graph $G=(V,E)$ is called $vertex$-$transitive$ if  $Aut(G)$
acts transitively on $V$. We say that $G$ is $edge$-$transitive$ if the group $Aut(G)$ acts transitively  on the edge set $E$, namely, for any $\{x, y\} ,   \{v, w\} \in E(G)$, there is some $a$ in $Aut(G)$,  such that $a(\{x, y\}) = \{v, w\}$.   We say that $G$ is $distance$-$transitive$ if  for all vertices $u, v, x, y$ of $G$ such that $d(u, v)=d(x, y)$, where $d(u, v)$ denotes the distance between the vertices $u$ and $v$  in $G$,  there is an automorphism $a$ in $Aut(G)$ such that  $a(u)=x$ and $a(v)=y.$

\section{Main Results}
 
Let $G=(V,E)$ be a graph with the vertex-set $V=\{ v_1,\dots,v_n \}$ and distance matrix  $D(G)=D=(d_{ij})_{n \times n}$, where $d_{ij}=d(v_i,v_j)$. Let $H \leq Aut(G)$  and 
 $ \pi = \{ w_1^H=C_1,\dots,w_m^H=C_m \} $ be the orbit partition of $H$, where $\{w_1,\dots,w_m  \} \subset V$. Let $Q=Q_{\pi}=(q_{ij})_{m \times m}$ be the matrix which its rows and columns is indexed by $\pi$ such that,
   $$q_{ij}=\sum _{w\in C_j} d(v,w), $$
     where $v$ is a fixed element in the cell $C_i$. It is easy to check that this sum is independent of $v$, that is, if $u \in C_i$, then $q_{ij}=\sum _{w\in C_j} d(v,w)=\sum _{w\in C_j} d(u,w)$. Hence, the matrix $Q$ is well defined. We call the matrix $Q$ the $quotient \ matrix $ of $ D $  over      $\pi  $.  We claim that every eigenvalue of $Q$ is an eigenvalue of the distance matrix $D$. In fact we have the following fact. 
 
  \begin{thm}Let $G=(V,E)$ be a graph with the distance matrix $D$. Let  $\pi$ be an orbit partition of $V$ with $m$ cells  and $Q$ be a quotient matrix of $D$ over $\pi$.  Then,  every eigenvalue of $Q$ is an eigenvalue of the distance matrix $D$.
  
   \end{thm}
   \begin{proof} 
     Let $\lambda$ be an eigenvalue of the quotient matrix $Q$ with a none zero eigenvector $f$. Let 
     $f(\pi_j)=x_j$. Thus for every $i, 1 \leq i \leq m, $ we have, \newline $\sum_{j=1}^{j=m}q_{ij} 
     f(\pi_j)$=$\sum_{j=1}^{j=m}q_{ij}x_j= \lambda x_i=\lambda f(\pi_i)$. \newline We define the function 
     $\hat{f}: V(G)\rightarrow \mathbb{R}$ by the rule $\hat{f}(v)=f(\pi_i)$ if and 
     only if $v \in \pi_i$. The function $\hat{f}$ is well defined since $\pi$ is a partition of $V(G)$. Also, $\hat{f}$ is a none zero function since the function $f$ is non zero.  If $v \in V(G)=V$, then there is a unique  $i$ such that $v \in \pi_i$. Now we have,\
 
 \

 $\sum_{w \in V} D_{vw}\hat{f}(w)=\sum_{w \in \pi_1} D_{vw}\hat{f}(w)+\sum_{w \in \pi_2} D_{vw}\hat{f}(w)+\dots +\sum_{w \in \pi_m} D_{vw}\hat{f}(w)=$\
 
 \

$\sum_{w \in \pi_1} D_{vw}f(\pi_1)+\sum_{w \in \pi_2} D_{vw}f(\pi_2)+\dots +\sum_{w \in \pi_m} D_{vw}f(\pi_m)=$\

\

$f(\pi_1) (\sum_{w \in \pi_1}D_{vw})+f(\pi_2)(\sum_{w \in \pi_2}D_{vw})+ \dots +f(\pi_m) (\sum_{w \in \pi_m}D_{vw})=$\

\

$f(\pi_1)q_{i1}+f(\pi_2)q_{i2}+\dots+f(\pi_m)q_{im}$=\

\

$\sum_{j=1}^{m} q_{ij} f(\pi_j)=\lambda f(\pi_i)=\lambda \hat{f}(v)$, since $ v\in \pi_i$.\\

Thus $\lambda$ is an eigenvalue of the matrix $D$. 

\end{proof}

By Theorem 2.1, we can find some of the distance eigenvalues of the graph $G$, but we can not  determine all the eigenvalues, since it is possible that $G$  has a  distance  eigenvalue 
$\theta$ such that  $\theta$ is not an eigenvalue of the  matrix $Q$. By the next theorem, we give a condition, holding that guaranties that the eigenvalue $\theta$ to be an eigenvalue  of the matrix  $Q$. 
 
\begin{prop} Let $G=(V,E)$ be a graph with the distance matrix $D$. Let  $\pi$ be an orbit partition of $V$ with $m$ cells  and $Q$ be the quotient matrix of $D$ over $\pi$.   Let $\theta$ be an eigenvalue of the distance matrix $D$, with the 
non zero eigenvector $f$ such that $f$ is constant on every cell of $\pi$.    Then  $\theta$ is an eigenvalue of the   matrix $Q$. 

\end{prop}

\begin{proof} 
We define the function $\tilde{f}: \pi \rightarrow \mathbb{R} $ by the rule $\tilde{f}(\pi_j)=f(v)$, where $v$ is an element in  the cell  $\pi_j$. Note that since $f$ is
 constant on the set $\pi_j$, then $\tilde{f}(\pi_j)$ is independent of $v \in \pi_j$,
  and hence $\tilde{f}$ is a well defined function. Also, since $f$ is non zero, then
   $\tilde{f}$ is non zero.
     If $u \in V$, since $f$ is an eigenvector with the eigenvalue $\theta $,  then we have,   
$$\sum_{w \in V} D_{uw}f(w)=\sum_{w \in V}d(u,w)f(w)=\theta f(u).$$

There is a unique $i$ such that $u \in \pi_i$.  Now we have,   
$$\sum_{w \in V}d(u,w)f(w)=\sum_{w \in \pi_1}d(u,w)f(w)+\sum_{w \in \pi_2}d(u,w)f(w)+\dots +\sum_{w \in \pi_m}d(u,w)f(w)$$
$$=\sum_{w \in \pi_1}d(u,w)\tilde{f}(\pi_1)+\sum_{w \in \pi_2}d(u,w)\tilde{f}(\pi_2)+\dots +\sum_{w \in \pi_m}d(u,w)\tilde{f}(\pi_m)$$
$$=\tilde{f}(\pi_1)(\sum_{w \in \pi_1}d(u,w))+\tilde{f}(\pi_2)(\sum_{w \in \pi_2}d(u,w))+\dots+\tilde{f}(\pi_m)(\sum_{w \in \pi_m}d(u,w))$$
$$=\tilde{f}(\pi_1)q_{i1}+\tilde{f}(\pi_2)q_{i2}+ \dots +\tilde{f}(\pi_m)q_{im}=\theta f(u)=\theta \tilde{f}(\pi_i). \ \ \ \ (*)$$
From ($*$) it follows that $\theta$ is an eigenvalue of the matrix $Q$ with the eigenvector $\tilde{f}.$
\end{proof}

Let $G=(V,E)$ be a graph with an adjacency matrix $A=(a_{vw})$ and automorphism group $\Gamma =Aut(G)$. We recall that every eigenvector $f$ of $G$ with the eigenvalue $\theta$ is a real function on $V$ such that $\sum_{w\in V}a_{vw}f(w)=\theta f(v)$, for every $v \in V$. If $g \in \Gamma$,  then the function $f^g$ defined by the rule $f^g(v)=f(v^g)$, $v \in V, $ is an eigenvector of $G$ with eigenvalue $\theta$ [7, chapt 9]. \

Let $D=(d_{vw})_{n \times n}$, $d_{vw}=d(v,w)$, be a distance matrix for $G$. Let $h$ be an eigenvector of $D$ with the eigenvalue $\lambda$.  Thus for every $v \in V$,  we have, 
  $$\sum_{w\in V}d_{vw}h(w)=\sum_{w\in V}d(v,w)h(w)=\lambda h(v).   \ 
   $$  
 We claim that if   $g \in Aut(\Gamma)$, then $h^g$ is an eigenvector of $D$ with the eigenvalue $\lambda$.  Note that
if $x,y \in V,$ then $d(x,y)=d(x^g,y^g)$.   We now have, 
$$\sum_{w \in V}d_{vw}h^g(w)=\sum_{w \in V}d(v,w)h(w^g)=\sum_{w \in V}d(v^g,w^g)h(w^g) =\lambda h(v^g)=\lambda h^g(v). $$
 The argument shows that $h^g$ is really an eigenvector of $D$ with the eigenvalue $\lambda$.
We now formally state the obtained result.

\begin{prop} Let $G=(V,E)$ be a   graph and $D$ be a distance matrix for $G$. Let $f$ be an eigenvector with the eigenvalue $\lambda$ for $D$. If $g$ is an automorphism of the graph $G$, then the function $f^g$ defined by the rule 
$f^g(v)=f(g^v)$, $v \in V$, is an eigenvector for $D$ with the eigenvalue $\lambda$.

\end{prop}
Let $G=(V,E)$ be a graph and $D$ be a distance matrix for $G$. Let $f$ be an eigenvector with the eigenvalue $\lambda$ for $D$. Let $H$ be a subgroup of $Aut(G). $ We can construct from $f$ an eigenvector $p$ with the eigenvalue $\lambda$ for $D$   such that $p$ is constant   on every orbit of $H$ on $V.$ In fact if we define the function $p$, by the rule, 
 $$p=\sum_{h \in H}f^h, \ \ \ \ \     (**)$$
  then from Proposition 2.3, 
it follows that if  $p\neq 0$, then  $p$ is an eigenvector of $D$ with the eigenvalue $\lambda$.  If $O=v^H,$ $v\in V$, is an orbit of $H$ on $V$, then $p$ is constant on $O$. Note that if $w \in O$, then $w=v^{h_1}$ for some $h_1 \in H$. Hence we have,   $$p(w)=\sum_{h \in H}f^h(w)=\sum_{h \in H}f(w^h)=\sum_{h \in H}f({v}^{h_1h})=\sum_{h \in H}f(v^h)=\sum_{h \in H}f^h(v)=p(v).$$
In other words, $p$ is an eigenvector of $D$ with the eigenvalue   $\lambda$ such 
that it is constant on every cell of $\pi$. Hence if $p \neq 0$, as we saw in proposition 2.2, we can construct the function $\tilde{p}$ from $p$ such that $\tilde{p}$ is an eigenvector for the matrix $Q$ with the eigenvalue $\lambda$,
  where $Q$ is a quotient matrix of $D$ over $\pi$.\\ 
We now can deduce that if $\lambda$ is an eigenvalue of the distance matrix $D$
with the eigenvector $f \neq 0$ such that $\lambda$ is not an eigenvalue of the matrix $Q$, then the function $p$ is the zero function. Thus we have the following result.

\begin{thm}Let $G=(V,E)$ be a vertex-transitive graph and $D$ be a distance matrix for $G$. Let $f \neq 0$ be an eigenvector with the eigenvalue $\lambda$ for $D$. Let $H$ be a subgroup of $Aut(G)$ and $\pi$ be its orbit partition on $V$ and $Q$ is a quotient matrix of $D$ over $\pi$. If $\lambda$ is not an eigenvalue of $Q$, then the sum of 
the values of $f$ on each cell of $\pi$ is zero.

\end{thm}

Now consider the function $p$ defined in ($**$). As we saw, it is  favorite that $p \neq 0$. In the next theorem we state some conditions, holding those 
guarantee that $p \neq 0$. 
  
\begin{prop} 
Let $G=(V,E)$ be a vertex-transitive graph and $D$ be a distance matrix for $G$.  Let $H$ be a subgroup of $Aut(G)$ with the orbit partition $\pi$ on $V$ such that $\pi$ has a singleton cell $\{x\}$. Let $\lambda$ be an eigenvalue of $D$. Then $\lambda $ is an eigenvalue  for the matrix $Q$, where $Q$ is a quotient matrix of $D$ over  $\pi$.

\end{prop}
\begin{proof}
Let $0 \neq f$ be an eigenvector with the eigenvalue $\lambda$ for $D$. Since $0 \neq f$, hence there is an element $w \in V$ such that $f(w) \neq 0.$ Since $G$ is a vertex-transitive graph, then $x^g=w,$ for some $g \in Aut(G).$ Let $t=f^g.$ Thus we have, 
$$t(x)=f^g(x)=f(x^g)=f(w)\neq 0.$$
 Let $p=\sum_{h \in H}t^h.$
 Then by Proposition 2.3,  $p$ is an eigenvector with the eigenvalue $\lambda$ for $D$ such that it  is constant on each cell of the partition $\pi$.   On the other hand we have,  $$p(x)=\sum_{h \in H}t^h(x)=\sum_{h \in H}t(x^h)=|H|t(x) \neq 0.$$
Hence  by proposition 2.2,  $\lambda$ is an eigenvalue for the matrix $Q.$ 
   
\end{proof}

From Theorem 2.1, Proposition 2.2, and Proposition 2.5, we obtain the following important result.

 \begin{thm} 
  Let $G=(V,E)$ be a vertex-transitive graph with the distance matrix $D$.
  Let $H$ be a subgroup of $Aut(G)$ with the orbit partition $\pi$ on $V$ such that $\pi$ has a singleton cell $\{x\}$. Let $Q=Q_{\pi}$ be a quotient   matrix of $D$ over $\pi$. Then the set of distinct eigenvalues of $D$ is equal to the set of distinct eigenvalues of $Q$. 
\end{thm} 

 In the sequel,  we will see   how Theorem  2.6, help us in finding the set of distance eigenvalues of the line graph of the crown graph.

Let $n \geq 3$ be an integer and $[n]= \{1,2,\dots,n \}$. Let $X= \{x_1,x_2,\dots,x_n  \}$ be an $n$-set disjoint from $[n]$.
 The crown graph $Cr(n)$ is a graph with the vertex-set $[n] \cup X$ and the edge-set $E_0=\{e_{ij}=\{i,x_j \} \ | \ i,j \in [n], i \neq j  \}.$   Thus, $L(Cr(n))$, the line graph of $Cr(n)$,  is a graph with the vertex-set $E_0$ in which two vertices $e_{ij}$ and $e_{rs} $ are adjacent if and only if $i=r$ or $j=s$.\
  
  Let $V=\{(i,j) \ | \ i,j \in [n], i \neq j  \}$. Let $G$ be a graph with the vertex-set $V$ in which two vertices $(i,j)$ and $(r,s)$ are adjacent if and only if $i=r$ or $j=s$.
  It is easy to check that the graph $G$ is isomorphic with the graph $L(Cr(n))$.
  Hence in the sequel we work on the graph $G$ and call it the line graph of  the crown graph and denote it by $L(Cr(n))$.\
  
It is easy to see that $L(Cr(3))$ is the cycle  graph $C_6$, which its structure is known and its line graph is again $C_6$. Hence in the rest of the paper we assume that $n\geq 4.$  It is easy to check that two non adjacent  vertices $(i,j)$ and $(r,s)$  are at distance 2 from each other whenever $i=s$ or $j=r$ or 
  $\{i,j\} \cap \{r,s\}=\varnothing$. Moreover, vertices  $(i,j)$ and $(j,i)$ are at distance 3 from each other ($P: (i,j),(x,j),(x,i),(j,i)$ is a shortest path between the vertices $(i,j)$ and $(j,i))$.   Thus the diameter of the graph $L(Cr(n))$ is 3. \
  
    Note that since the crown graph $Cr(n)$ is a regular graph and its adjacency spectrum is known 
 [4], hence the adjacency spectrum of its line graph, that is, $L(Cr(n))$ is known [3,5,7].
  \begin{rem} Although the crown graph $Cr(n)$ is a distance-transitive graph (and consequently it is distance-regular [4,9]), it is easy to check that the graph $L(Cr(n))$ is not distance-regular. Hence we can not use the theory of distance -regular graphs for determining the set of distance eigenvalues of this graph.
  \end{rem} 
   Figure 1, shows the graph $L(Cr(4))$. Note that in this figure the vertex $(i,j)$
 is denoted by $ij$. \
 
\definecolor{qqqqff}{rgb}{0.,0.,1.}
\begin{tikzpicture}[line cap=round,line join=round,>=triangle 45,x=0.6684255715985469cm,y=.75cm]
\clip(3.371,-7.1306) rectangle (22.0292,3.0576);
\draw (9.8566,2.6462) node[anchor=north west] {41};
\draw (14.7692,2.5736) node[anchor=north west] {42};
\draw (16.8746,-0.0642) node[anchor=north west] {12};
\draw (16.8746,-2.5568) node[anchor=north west] {13};
\draw (15.8098,-5.2188) node[anchor=north west] {14};
\draw (9.9292,-5.4366) node[anchor=north west] {34};
\draw (7.5092,-2.7746) node[anchor=north west] {32};
\draw (6.9914,-0.3788) node[anchor=north west] {31};
\draw (9.1114,0.3714) node[anchor=north west] {21};
\draw (12.8332,-4.5896) node[anchor=north west] {24};
\draw (14.8902,-2.0002) node[anchor=north west] {23};
\draw (10.42,1.8)-- (14.7692,1.8476);
\draw (14.7692,1.8476)-- (16.52,-0.7);
\draw (16.52,-0.7)-- (16.5,-2.94);
\draw (16.5,-2.94)-- (15.08,-5.24);
\draw (15.08,-5.24)-- (10.54,-5.22);
\draw (10.54,-5.22)-- (8.22,-2.8);
\draw (8.22,-2.8)-- (8.16,-0.7);
\draw (8.16,-0.7)-- (10.42,1.8);
\draw (10.42,1.8)-- (10.32,-0.12);
\draw (10.32,-0.12)-- (8.16,-0.7);
\draw (12.86,-4.34)-- (10.54,-5.22);
\draw (12.86,-4.34)-- (15.08,-5.24);
\draw (14.04,-0.66)-- (14.48,-2.46);
\draw (14.7692,1.8476)-- (14.04,-0.66);
\draw (14.2318,0.1326) node[anchor=north west] {43};
\draw (16.5,-2.94)-- (14.04,-0.66);
\draw (14.48,-2.46)-- (16.5,-2.94);
\draw (10.32,-0.12)-- (14.48,-2.46);
\draw (16.52,-0.7)-- (15.08,-5.24);
\draw (10.32,-0.12)-- (12.86,-4.34);
\draw (16.52,-0.7)-- (8.22,-2.8);
\draw (14.48,-2.46)-- (12.86,-4.34);
\draw (10.42,1.8)-- (14.04,-0.66);
\draw (14.7692,1.8476)-- (8.22,-2.8);
\draw (8.16,-0.7)-- (10.54,-5.22);
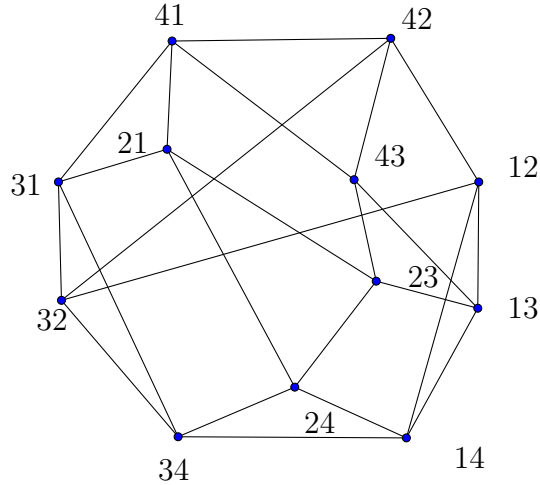
\draw (8.4328,-6.078) node[anchor=north west] {Figure 1. The graph $L(Cr(4))$};
\begin{scriptsize}
\draw [fill=qqqqff] (10.42,1.8) circle (1.5pt);
\draw [fill=qqqqff] (14.7692,1.8476) circle (1.5pt);
\draw [fill=qqqqff] (16.52,-0.7) circle (1.5pt);
\draw [fill=qqqqff] (16.5,-2.94) circle (1.5pt);
\draw [fill=qqqqff] (15.08,-5.24) circle (1.5pt);
\draw [fill=qqqqff] (10.54,-5.22) circle (1.5pt);
\draw [fill=qqqqff] (8.22,-2.8) circle (1.5pt);
\draw [fill=qqqqff] (8.16,-0.7) circle (1.5pt);
\draw [fill=qqqqff] (10.32,-0.12) circle (1.5pt);
\draw [fill=qqqqff] (12.86,-4.34) circle (1.5pt);
\draw [fill=qqqqff] (14.04,-0.66) circle (1.5pt);
\draw [fill=qqqqff] (14.48,-2.46) circle (1.5pt);
\end{scriptsize}
\end{tikzpicture}
\

Since the graph $Cr(n)$ is distance-transitive, hence it is edge-transitive. Thus the graph $L(Cr(n))$ is a vertex-transitive graph. For each $\alpha \in Sym([n])$, let $f_{\alpha}$ be the function on the vertex-set of the graph $L(Cr(n))$ defined by the rule,  $f_{\alpha}(i,j)=(\alpha(i),\alpha(j))$. Let $\beta$ be the function on the vertex-set of the graph $L(Cr(n))$ defined by the rule,  $\beta(i,j)=(j,i)$. Now we can check that $Aut(L(Cr(n))) \cong Sym([n]) \times \langle \beta \rangle$ [9,11,13], where $ \langle \beta \rangle$ is the subgroup generated by the automorphism $\beta$ in the automorphism group of the graph $L(Cr(n))$. 
\begin{thm}
Let $n>3$ be an integer. Then the line graph of the crown graph, that is, the graph $L(Cr(n))$ is a distance  integral graph with distinct  distance eigenvalues,  $-n-1, -n+3, - 1, 1,  2n^2-4n+3$.
\end{thm} 

\begin{proof}
Let $D$ be a distance matrix of the graph $L(Cr(n))$.
In the first step, we  proceed to construct an orbit partition $\pi$  for  the vertex-set of $L(Cr(n))$,   such that this partition has a singleton cell. 
If we construct such a partition, then   since the graph $L(Cr(n))$  is a vertex-transitive graph, then  by Theorem  2.6,  every   distance eigenvalue of the graph $L(Cr(n))$ is an eigenvalue of the matrix $Q$ and vice versa, where $Q$ is the quotient matrix of $D$ over $\pi$.   
Let $ H = \{  f_{ \alpha } \  | \    \alpha \in Sym([n]),\  \alpha(1) = 1, \   \alpha(2) = 2 \} $.  Then $H$ is a subgroup of $ Aut( L(Cr(n)))$, the automorphism group of $L(Cr(n)) $.  Let $ H_1 = \{ \alpha \ | \  f_{ \alpha}     \in H  \}   $.  Note that $H_1$ is a subgroup of $Sym([n])$ isomorphic with $Sym([n-2])$. In the sequel,  we want to determine the orbit partition of  the subgroup $ H $. In fact,  $H$ generates the following orbits; \\\\
$O_1 =H((1,2)) = \{   h((1,2))  
 \ | \ h \in  H \}$ = $ \{ ( \alpha(1), \alpha(2) )  
 \ | \ \alpha \in H_1\} $= $\{(1,2)\}$. \\\\
$O_2 =H((1,3)) = \{   h((1,3))  
 \ |  \  h \in  
H \}$ = $ \{ ( \alpha(1), \alpha(3) )  
 \ |  \  \alpha \in H_1\} $= $\{  (1,i)  \ | \  3 \leq i \leq n \}$. \\\\
 $O_3 =H((3,1)) = \{   h((3,1))  
 \ | \   h \in 
H \}$ = $ \{ ( \alpha(3), \alpha(1))  
 \ |  \  \alpha \in H_1\} $= $\{  (i,1) \  | \  3 \leq i \leq n \}$. \\\\
$O_4 =H((2,1)) = \{   h((2,1))  
 \ | \  h\in  
H \}$ = $ \{ ( \alpha(2), \alpha(1))  
 \ | \  \alpha \in H_1\} $= $\{ (2,1) \}$. \\\\
 $O_5 =H((2,3)) = \{   h((2,3))  
 \ |  \  h \in  
H \}$ = $ \{ ( \alpha(2), \alpha(3) ) 
 \ | \  \alpha \in H_1\} $= $\{  (2,i) \  |  \ 3 \leq i \leq n \}$. \\\\
$O_6 =H((3,2)) = \{   h((3,2))  
 \ | \  h \in  
H \}$ = $ \{ ( \alpha(3), \alpha(2) ) 
\  | \  \alpha \in H_1\} $= $\{  (i,2)  \ |  \ 3 \leq i \leq n \}$.\\\\
$O_7 =H((3,4)) = \{   h((3,4))  
\ | \  h \in 
H \}$ = $ \{ ( \alpha(3),\alpha(4)) 
\ | \  \alpha \in H_1\} $= $\{  (i,j) \  | \  3 \leq i,j \leq n, \ i \neq j \}$. \\\\
If we let $ \pi = \{ O_1 , O_2, O_3, \dots ,O_7 \} $, then  $ O_1 \cup O_2 \cup \dots \cup O_7 = V= V(L(Cr(n))) $.
Let $Q=(q_{ij})_{7 \times 7}$, be the quotient matrix of $D$ over $\pi$, that is,      $q_{ij}$ is the sum of the distances of  a vertex in the cell $O_i$ from all the vertices in the cell $O_j$.  Then  the following hold. \\\\
$q_{11}=0$. \\\\
$q_{12}= n-2$.    Because  the vertex $(1,2) \in O_1$ is adjacent to all the $n-2$ vertices in $O_2$.\\\\
$q_{13}= 2(n-2)=2n-4$.    Because  the distance of the vertex $(1,2) \in O_1$ is 2 from every vertex in $O_3$. Note that $|O_3|=n-2$. \\\\
$q_{14}=3$.    Because    the distance of the vertex $(1,2) \in O_1$ is 3 from the  vertex  $(2,1) \in O_4$.  \\\\
$q_{15}=2(n-2)= 2n-4$.   Because   the distance of the vertex $(1,2)\in O_1$ is 2 from every vertex in $O_5$. Note that $|O_5|=n-2$.\\\\
$q_{16}=n-2$.   Because   the distance of the vertex $(1,2) \in O_1$ is 1 from every vertex in $O_6$.\\\\
$q_{17}=2(n-2)(n-3)$.  Because  the distance of the vertex $(1,2) \in O_1$ is 2 from every vertex in $O_7$. Note that $|O_7|=(n-2)(n-3). $\\\\
$q_{21}=1.$  Because   the distance of the vertex $(1,3) \in O_2$ is 1 from the vertex $(1,2)$ in  $O_1$. \\\\
$q_{22}=n-3$.  Because   the distance of the vertex $(1,3) \in O_2$ is 1 from each of  the other vertices in $O_2$.  Note that $|O_2|=n-2.$ \\\\
$q_{23}=2n-3$.   Because   the distance of the vertex $(1,3) \in O_2$ is 3 from the vertex $(3,1)$, and is 2  from each of  the other vertices in $O_3$. Hence we have $q_{23}=3+2(n-3)=2n-3$.  Note that $|O_3|=n-2.$ \\\\
$q_{24}$=2.  Because  the distance of the vertex $(1,3) \in O_2$ is 2 from the vertex $(2,1)$ in $O_4$. \\\\
$q_{25}=2n-5.$  Because   the distance of the vertex $(1,3) \in O_2$ is 1 from the vertex $(2,3),$ and is 2  from each of the other vertices in  $O_5.$ Hence we have $q_{25}=1+2(n-3)=2n-5$. Note that $|O_5|=n-2$. \\\\
$q_{26}=2n-4.$  Because  the distance of the vertex $(1,3) \in O_2$ is 2 from every vertex in $O_6.$ Note that $|O_6|=n-2$. \\\\
$q_{27}=(n-3)(2n-5)$.  Because   the distance of the vertex $(1,3) \in O_2$ from every vertex of the form $(j,3)$, $4 \leq j \leq n$,  in  $O_7$ is 1 and from each of  the other vertices in  $O_7$ is 2.  Hence  we have $q_{27}=(n-3)+2((n-3)+(n-4)(n-3))=(n-3)+2(n-3)(n-3)$=$(n-3)(2n-5)$.  \\\\
$q_{31}=2.$  Because   the distance of the vertex $(3,1) \in O_3$ is 2 from the vertex $(1,2)$   in $O_1$. \\\\
$q_{32}=2n-3.$  Because   the distance of the vertex $(3,1) \in O_3$ is 3 from the vertex $(1,3)$ in $O_2$  and from each of the other vertices is 2. Hence we have $q_{32}=3+2(n-3)$. \\\\
$q_{33}=n-3.$  Because $|O_3|=n-2$ and the  vertex $(3,1) \in O_3$ is adjacent to each of the other vertices in $O_3.$ \\\\
$q_{34}=1.$\\\\ 
$q_{35}=2n-4.$  Because   the distance of the vertex $(3,1) \in O_3$ is 2 from every vertex in $O_5$. \\\\
$q_{36}=2n-5.$  Because  the distance of the vertex $(3,1) \in O_3$ is 1  from the vertex $(3,2)$ in $O_6$  and from each of the other vertices is 2. Hence we have $q_{36}=1+2(n-3)=2n-5$. \\\\
$q_{37}=(n-3)(2n-5).$  Because   the distance of the vertex $(3,1) \in O_3$ from every vertex of the form $(3,j)$, $4 \leq j \leq n$,  in  $O_7$ is 1 and from each of  the other vertices in  $O_7$ is 2.  Hence  we have $q_{37}=(n-3)+2((n-3)(n-3))=$$(n-3)(2n-5)$. \\\\
$q_{41}=3.$  \\\\
$q_{42}=2n-4.$  Because   the distance of the vertex $(2,1) \in O_4$ is 2 from every vertex in $O_2$, hence we have $q_{42}=2(n-2)$. \\\\
$q_{43}=n-2.$  Because  the vertex  $(2,1) \in O_4$ is adjacent to every vertex in $O_3.$ \\\\
$q_{44}=0.$  \\\\  
   $q_{45}=n-2.$ Because  the vertex  $(2,1) \in O_4$ is adjacent to every vertex in $O_5.$ \\\\
$q_{46}=2n-4.$  Because  the distance of the vertex $(1,2) \in O_4$ is 2 from every vertex in $O_6$. Hence we have $q_{46}=2(n-2). $ \\\\    
 $q_{47}=2(n-2)(n-3).$  Because   the distance of the vertex $(2,1) \in O_4$ is 2 from every vertex in $O_7$.  Since $|O_7|=(n-2)(n-3)$, thus we have $q_{47}=2(n-2)(n-3).$ \\\\    
$q_{51}=2.$  Because  the distance of the vertex $(2,3) \in O_5$ is 2 from the vertex $(1,2)$   in $O_1$. \\\\  
 $q_{52}=2n-5.$  Because  the distance of the vertex $(2,3) \in O_5$ is 1 from the vertex $(1,3)$ in $O_2$, and is 2 from each of the other vertices in $O_2$. Hence we have $q_{52}=1+2(n-3)=2n-5.$ \\\\    
 $q_{53}=2n-4.$   Because   the distance of the vertex $(2,3) \in O_5$ is 2 from every vertex in $O_3$. \\\\            
 $q_{54}=1.$  Because the vertex  $(2,3) \in O_5$ is adjacent to vertex $(2,1)$ in $O_4$. \\\\                 
$q_{55}=n-3.$  Because    the vertex $(2,3) \in O_5$ is adjacent to every other vertex in $O_5$. \\\\ 
$q_{56}=2n-3.$  Because  the distance of the vertex $(2,3) \in O_5$ is 3 from the 
vertex $(3,2)$ in $O_6$ and is 2 from  every other  vertex in $O_6$. Thus we have
$q_{56}=3+2(n-3).$ \\\\ 
 $q_{57}=(n-3)(2n-5).$  Because   the distance of the vertex $(2,3) \in O_5$ from every vertex of the form $(j,3)$, $4 \leq j \leq n$,  in  $O_7$ is 1 and from each of  the other vertices in  $O_7$ is 2.  Hence  we have $q_{57}=(n-3)+2((n-3)+(n-4)(n-3))=(n-3)+2(n-3)(n-3)$=$(n-3)(2n-5)$. \\\\                           
$q_{61}=1.$  Because  the   vertex $(3,2) \in O_6$ is  adjacent to the vertex $(1,2)$ in $O_1$. \\\\                                                       
 $q_{62}=2n-4.$   Because   the distance of the vertex $(3,2) \in O_6$ is 2 from every vertex in $O_2$. Hence we have $q_{62}=2(n-2)$.\\\\                                                                                 
 $q_{63}=2n-5.$  Because  the distance of the vertex $(3,2) \in O_6$ is 1 from the 
vertex $(3,1)$ in $O_3$ and is 2 from  every other  vertex in $O_3$. Thus we have
$q_{63}=1+2(n-3).$\\\\
 $q_{64}=2.$  Because  the distance of the vertex $(3,2) \in O_6$ is 2 from the vertex  $(2,1)$ in $O_4$. \\\\
  $q_{65}=2n-3.$  Because  the distance of the vertex $(3,2) \in O_6$ is 3 from the vertex $(2,3)$ in $O_5$, and is 2 from every  other vertex in $O_5$.  Thus we have,  $q_{65}=3+2(n-3). $\\\\
   $q_{66}=n-3.$  Because the vertex $(3,2) \in O_6$  is adjacent to every other vertex in $O_6$. Note that $|O_6|=n-2.$\\\\
    $q_{67}=(n-3)(2n-5)$.   Because  the distance of the vertex $(3,2) \in O_6$ from every vertex of the form $(3,j)$, $4 \leq j \leq n$,  in  $O_7$ is 1 and from each of  the other vertices in  $O_7$ is 2.  Hence  we have $q_{67}=(n-3)+2((n-3)+(n-4)(n-3))=(n-3)+2(n-3)(n-3)$=$(n-3)(2n-5)$.\\\\
     $q_{71}=2.$  Because   the distance of the vertex $(3,4) \in O_7$ is 2 from the vertex  $(1,2)$ in $O_1$. \\\\
     $q_{72}=2n-5.$  Because  the distance of the vertex $(3,4) \in O_7$ is 1 from the vertex $(2,4)$ and is 2 from every  other vertex in $O_2$. Thus we have $q_{72}=1+2(n-3)$.\\\\
      $q_{73}=2n-5.$  Because the distance of the vertex $(3,4) \in O_7$ is 1 from the vertex $(3,1), $ and is 2 from every  other vertex in $O_3$. Thus we have $q_{73}=1+2(n-3)$. \\\\
$q_{74}=2.$  Because  the distance of the vertex $(3,4) \in O_7$ is 2 from the vertex  $(2,1)$ in $O_4$.\\\\
 $q_{75}=2n-5.$  Because the distance of the vertex $(3,4) \in O_7$ is 1 from the 
vertex $(2,4)$ in $O_5$ and is 2 from  every other  vertex in $O_5$. Thus we have
$q_{75}=1+2(n-3).$  \\\\
 $q_{76}=2n-5.$  Because  the distance of the vertex $(3,4) \in O_7$ is 1 from the 
vertex $(3,2)$ in $O_6, $ and is 2 from  every other  vertex in $O_6$. Thus we have
$q_{76}=1+2(n-3).$ \\\\
$q_{77}=2(n-4)(n-2)+3.$  Because  the distance of the vertex $(3,4) \in O_7$  from every vertex in $O_7$ of the form $(3,j),$ $5 \leq j \leq n$, is 1, and  from every vertex in $O_7$ of the form $(i,4),$ $5 \leq i \leq n$, is 1, and from the vertex $(4,3)$ is 3, and from every vertex in $O_7$ of the form $(i,j),$ $4 \leq i \leq n$, $j \neq 4$, is 2. Hence we have $q_{77}=(n-4)+(n-4)+3+2(n-4)(n-3)=2(n-4)(1+n-3)+3$=$2(n-4)(n-2)+3$. \\\\
Therefore,  we obtain the following quotient matrix $ Q$ of $D$ over $\pi$.  
$$Q=\begin{pmatrix}
0 & n-2 & 2n-4 & 3 & 2n-4 & n-2 & 2(n-2)(n-3) \\
 1 & n-3 & 2n-3 & 2 & 2n-5 & 2n-4 & (n-3)(2n-5) \\
 2 & 2n-3 & n-3 & 1 & 2n-4 & 2n-5 & (n-3)(2n-5) \\
 3 & 2n-4 & n-2 & 0 & n-2 & 2n-4 & 2(n-2)(n-3) \\
 2 & 2n-5 & 2n-4 & 1 & n-3 & 2n-3 & (n-3)(2n-5) \\
 1 & 2n-4 & 2n-5 & 2 & 2n-3 & n-3 & (n-3)(2n-5) \\
 2 & 2n-5 & 2n-5 & 2 & 2n-5 & 2n-5 & 2(n-4)(n-2)+3
\end{pmatrix}$$ \\
We can use the Wolfram Mathematica [15] for finding the eigenvalues of the matrix $Q$. Using  the Wolfram Mathematica  we have,  \newline \newline
$Q= \{\{0,n-2,2n-4,3,2n-4,n-2,2(n-2)(n-3)\},\\
\{1,n-3,2n-3,2,2n-5,2n-4,(n-3)(2n-5)\},\\
\{2,2n-3,n-3,1,2n-4,2n-5,(n-3)(2n-5)\},\\
\{3,2n-4,n-2,0,n-2,2n-4,2(n-2)(n-3)\},\\
\{2,2n-5,2n-4,1,n-3,2n-3,(n-3)(2n-5)\},\\
\{1,2n-4,2n-5,2,2n-3,n-3,(n-3)(2n-5)\},\\
\{2,2n-5,2n-5,2,2n-5,2n-5,2(n-4)(n-2)+3\}\}\\\\
{\bf Eigenvalues[Q]}=  \{ -1,1,-1-n,-1-n,3-n,3-n,3-4 n+2 n^2 \}. $\\\\
  We now conclude that the set  $  \{ -n-1,-n+3,-1,1, 2n^2-4n+3 \},  $ is the set of all  distinct eigenvalues of the graph $ L(Cr(n)) $.
\end{proof} 
  
\begin{rem} It is not true that if $G=(V,E)$ is an integral distance-transitive graph, then its line graph $L(G)$ is distance integral. Using Wolfrom Mathematica [15], one can see that the Johnson graph $J(6,2)$ which is an integral distance-transitive graph, is distance integral.  But its line graph, that is,   $L(J(6,2))$ is not distance integral.

\end{rem}

\end{document}